\documentclass[psamsfonts, reqno]{amsart}

\usepackage[english]{babel}
\usepackage{epsfig,epic,eepic}
\usepackage{amssymb,algorithmic}
\usepackage[hyphens]{url}

\theoremstyle{plain}
\newtheorem{theorem}{Theorem}[section]
\newtheorem{lemma}[theorem]{Lemma}

\newtheorem{theorem*}{Theorem}

\theoremstyle{definition}
\newtheorem{definition}[theorem]{Definition}
\newtheorem{conjecture*}{Conjecture}

\def\B{{\mathcal B}}

\def\R{{\mathbb R}}

\def\A{{\mathbf A}}
\def\BB{{\mathbf B}}
\def\cm{{\alpha}}
\def\cM{{\beta}}

\def\Pr{{\bf Pr}}
\def\E{{\bf E}}
\def\Var{{\bf Var}}

\ifx\hyperlink\undefined
\newcommand{\hrlb}{ }
\else
\newcommand{\hrlb}{\\}
\fi

\begin{document}

\author{Gady Kozma}
\address{GK: Institute for Advanced Study, 1 Einstein Dr., Princeton NJ
  08540, USA}
\email{gady@ias.edu}
\author{Zvi Lotker}
\address{ZL: Office Centrum voor Wiskunde en Informatica 
Kruislaan
413, NL-1098 SJ Amsterdam,
The Netherlands.}
\email{lotker@cwi.nl}
\author{Gideon Stupp}
\address{GS: MIT Computer Science and Artificial Intelligence Laboratory,
The Stata Center, 32 Vassar St., Cambridge, MA 02139, USA.}
\email{gstupp@theory.csail.mit.edu}

\title{On the connectivity of the Poisson process on fractals}


\begin{abstract}

  For a measure $\mu$ supported on a compact connected subset of a
  Euclidean space which satisfies a
  uniform $d$-dimensional decay of the volume of balls of the type
  \begin{equation}
  \cm\delta^d \leq \mu(B(x,\delta)) \leq \cM\delta^d\label{uniform}
  \end{equation}
  we show that the maximal edge in the minimum spanning tree of $n$
  indepndent samples from $\mu$ is, with high probability
  $$\approx\left(\frac{\log n}{n}\right)^{1/d}.$$

\end{abstract}

\maketitle
\section{Introduction}
The laws governing the behavior of the minimum
spanning tree (MST) on points in the unit disk and ball
have been thoroughly
researched. See \cite{S90} for a survey on the worst case. For the
average case (or points taken randomly), extremely fine results are
known. See \cite{A96} for the central limit theorem for the total
length, \cite{ABNW99, AB99} for the dimension of a typical path,
\cite[chapter 6]{AS02} for an ``objective'' approach and \cite{W} for
intriguing simulation results. See also the book \cite{MR96} for the
strongly related continuum percolation.

The same questions on other metric spaces, in particular on fractals,
have seen much less research. In \cite{KLS} we studied the worst case
problem for the total weighted length. Here we switch to the average
case and are interested in the length of the longest edge, which, by
the greedy algorithm, is the same as the connectivity threshold
i.e.~the minimal number $r$ such that the graph in which two points are
connected if and only if their metric distance is $\leq r$, is
connected. In the setting of a ball in $\R^d$ this is known to be, with high
probability
$$\approx \left(\frac{\log n}{n}\right)^{1/d}$$
where $\approx$ means that the ratio of the two quantities is bounded
between two absolute constants. We wish to extend this result to
fractal sets.

\begin{figure}[h]
\begin{center}
\psfig{figure=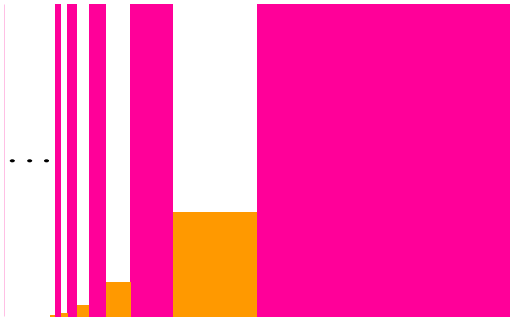,scale=0.5}\caption{ The set F (not to scale).\label{set}}
\end{center}
\end{figure}

Clearly, to get any kind of estimate one has to assume the fractal is
connected. Further, it is clear that some kind of regularity is
needed. To see why, it might be instructive to consider the following
example: in $\R^2$ take the set
 $F=\cup_{i=1}^\infty(A_i \cup
B_i)$, consisting of a set of ``thick'' vertical slabs
$A_i=\left[\frac{1}{2i},\frac{1}{(2i-1)}\right]\times[0,1]$
connected by ``thin'' horizontal bridges
$B_i=\left[\frac{1}{(2i+1)},\frac{1}{2i}\right]\times\left[0,2^{-i}\right]$;
see figure~\ref{set}. Take the normalized Lebesgue measure on this
set. Now take $n$ random points and connect them by their MST. Since
the bridges become thin very fast, there will be no points in the
very thin bridges that connect the slabs starting from
$[1/2\lfloor\log n\rfloor, 1/(2\lfloor\log
  n\rfloor-1)]\times[0,1]$. Thus, the MST will contain edges that 
are $\geq c/\log n$ long. This is despite the fact that this set has
Hausdorff dimension $2$ and in fact is a monofractal (a degenerate
multifractal spectrum) which indicates strong regularity in some
sense.

\begin{definition}
  The metric probability space $(F,\rho,\mu)$ is \emph{semi uniform of
  dimension $d$} if
  there exists some numbers $\delta_0$ and $0<\cm\leq\cM<\infty$ such
  that for every point $x\in F$ and for every $0<\delta<\delta_0$,
$$\cm\delta^d\leq\mu(B(x,\delta))\leq\cM\delta^d$$
where $B(x,r)$ is the ball of radius $r$ centered at $x$ in the metric $\rho$.
We call $\delta_0$, $\cm$ and $\cM$ the \emph{parameters} of the
metric probability space.
\end{definition}

Examples for such a semi
uniform space are the Cantor dust or the Sierpinski carpet, endowed
with the metric of their embedding in $\R^n$ and the natural
probability measures~\cite{F90}. We show that for such probability
spaces the classical asymptotic laws of maximal MST edge still hold
with fractional powers, if the space is connected. For example, in the
case of Sierpinski's carpet, with dimension $d=\log 8/\log 3$, the
longest edge of the MST is $\approx\left(\frac{\log
  n}{n}\right)^{1/d}$ with probability $1-\epsilon$ for all
$\epsilon>0$. Formally, the statement is

\begin{theorem} Let $F$ be a compact connected subset of $\R^k$ and
  let $\mu$ be a semi-uniform measure of dimension $d$ on $F$.
Then there exist two constants $C>c>0$ such that for all
$\epsilon>0$ and for any $m$ sufficiently large (depending on
$\epsilon$), if $X_1,\dotsc,X_m$ are independent samples from
$\mu$ then

  \begin{enumerate}
  \item If $r>C(\frac{\log m}{m})^{1/d}$ then the set
    $\bigcup_i B(X_i,r)$ is connected with probability $1-\epsilon$.
  \item if $r<c(\frac{\log m}{m})^{1/d}$ then the set $\bigcup_i
    B(X_i,r)$ is connected with probability $\epsilon$.
  \end{enumerate}
\end{theorem}

Here and below we shall use the convention that $c$ and $C$ denote constants
whose value might change from formula to formula and even inside the
same formula. $c$ and $C$ might depend on $d$, $\cm$ and $\cM$. $c$
will usually denote constants which are ``small enough'' and $C$
constants which are ``large enough''.
\section{Proof}

To prove the first part of the theorem choose
$\delta:=\left(2\log m/(\cm m)\right)^{1/d}$ and let
$\{B(p_i,\delta)\}_{i=1}^{N(\delta)}$ be a maximal set of disjoint
balls with centers $p_i\in F$.  Then for each $X_i$,
$$
\Pr[X_i\in B(p_j,\delta)]\geq 2\frac{\log m}{m}\enspace .
$$
Here and below $\Pr$ stands for probability with respect to the product
measure $\prod_{i=1}^{N(\delta)}\mu$.
Denote by $\A_\delta$ the event that there is at least one of $X_i$ in each
of $B(p_j,\delta)$.

\begin{lemma}
  For large enough $m$, event $\A_\delta$ occurs with probability $1-\epsilon$.
\end{lemma}

\begin{proof}
  Had all the measures of all the balls been equal, we could have
  immediatly used coupon collector. However, the measures of
  the balls are only equal up to a constant factor. Still, for all $p_j$,
  $$2\frac{\log m}{m}\leq \Pr[B(p_j,\delta)]\enspace .$$
  Take $m/(2\log m)$ equiprobable
  coupons and denote by $\A'$ the probability of picking at least one
  coupon out of $\A'$ after $m$ trials. Then by coupon collector
  \[\Pr[\A_\delta]\geq\Pr[\A']\geq 1-\epsilon\enspace .\qedhere\]
\end{proof}

\begin{lemma}
  Conditioned on $\A_\delta$, $\bigcup_i B(X_i,3\delta)$ is connected.
\end{lemma}
\begin{proof}%
As is well known (and easy to see), the maximality of the family
$B(p_i,\delta)$ implies that the family $B(p_i,2\delta)$ is
a cover of $F$. Since we are given that the event $\A_\delta$
happens it follows that the set $\bigcup B(X_i,3\delta)$ contains the
set $\bigcup B(p_i,2\delta)$ and therefore it also contains
$F$. Now the lemma follows directly from the definition of
connectivity. Suppose to the contradiction that $\bigcup
B(X_i,3\delta)$ is separated i.e.~there exist two open sets $U,V$
s.t.~$U\cup V = \bigcup B(X_i,3\delta)$ and $U\cap V=\emptyset$. Using
$U,V$ it follows that $F$ is also separated, and therefore a
contradiction.
\end{proof}
\subsection{}
With the first part of the theorem demonstrated, let us proceed to the
    second. This time choose $\delta=\left(\frac{\log m}{2\cM
    m}\right)^{1/d}$ and have
$\{B(p_i,\delta)\}_{i=1}^{N(\delta)}$ be a maximal set of
disjoint balls with centers $p_i\in F$.  Then for each $X_i$,
\[
\frac{\cm}{\cM}\cdot\frac{\log m}{2m} \leq
\Pr[X_i\in B(p_j,\delta)]\leq \frac{\log m}{2m}\enspace .
\]
Denote by $\BB_{\delta}$ the event that there are at least $\log m$
balls in $\{B(p_j,\delta)\}$ that contain exactly a single point $X_i$.  
Denote by $Y_j$ the indicator to the event that the ball
$B(p_j,\delta)$ contains exactly one point $X_i$ and have
$Y=\sum Y_j$ be the random variable that counts the number of such
balls in $\{B(p_j,\delta)\}$. It is not difficult to see that
$$\E(Y_i)\geq \frac{\cm}{\cM}\cdot\frac{\log m}{2}
     \left(1-\frac{\log m}{2m}\right)^{m-1}
  = \frac{\log m}{\sqrt{m}}\left(\frac{\cm}{2\cM}+o(1)\right)
\enspace ,$$
so for $m$ sufficiently large,
\begin{equation}
\E(Y) \geq c\sqrt{m}\log m\enspace .\label{EYbig}
\end{equation}

\begin{lemma}\label{lem:ugly}
For $m$ sufficiently large,
\[
\Var(Y)<C\frac{\E(Y)^2}{\sqrt{m}}\enspace .
\]
\end{lemma}

\begin{proof}





For any $1\leq i\neq j\leq N(\delta)$ examine the variables $Y_i$ and
$Y_j$. Denote $q_i = \mu(B(p_i,\delta))$. Then, because $q_i \leq
C\log m/m$ we have
\begin{align*}
\E(Y_i)&=mq_i(1-q_i)^{m-1}=mq_i\exp\left((-q_i+O(q_i^2))(m-1)\right)= \\
&=mq_i\exp\left(-q_im+O\left(\frac{\log^2 m}{m}\right)\right)=\\
&=mq_i\exp(-mq_i)\left(1+O\left(\frac{\log^2 m}{m}\right)\right)
\enspace .
\end{align*}
A similar calculation gives
\begin{align*}
\E(Y_i Y_j)&=m(m-1)q_i q_j(1-q_i-q_j)^{m-2}=\\
&=m^2 q_i q_j \exp(-mq_i-mq_j)
  \left(1+O\left(\frac{\log^2 m}{m}\right)\right)\enspace .
\end{align*}
Combining both we get, for $m$ sufficiently large,
\begin{align*}
\Var(Y)&=\E(Y^2)-\E(Y)^2=
  \sum_{i=1}^{N(\delta)}\sum_{j=1}^{N(\delta)}
  \E(Y_i Y_j) - \E(Y_i)\E(Y_j)\leq \\
&\leq \sum_{i=1}^{N(\delta)}\E(Y_i) + \sum_{i\neq j} m^2q_i q_j
  \exp(-mq_i -mq_j) \cdot C\frac{\log^2 m}{m} \leq\\
&\leq \E(Y)+C\frac{\log^2 m}{m}\E(Y)^2\enspace .
\end{align*}
This ends the lemma by virtue of (\ref{EYbig}).
\end{proof}

\begin{lemma}
  For large enough $m$, event $\BB_\delta$ occurs with probability $1-\epsilon$.
\end{lemma}

\begin{proof}
  By Lemma~\ref{lem:ugly} and the Chebyshev inequality it follows that
\[
  \Pr\left[\left|Y-\E(Y)\right|>\frac{1}{2}\E(Y)\right]\leq \frac
  {4}{\sqrt{m}} \enspace.\qedhere
\]
\end{proof}

\begin{lemma}
  Conditioned on $\BB_\delta$, $\bigcup_i\B(X_i,\frac{1}{2}\min_j(\delta))$ is
  disconnected with probability $1-\epsilon$.
\end{lemma}
\begin{proof}
  Given that $B(p_j,\delta)$ is a ball containing just a single point
  $X_i$,  it follows that there is a
  constant probability that this point will be in
  $B(p_j,\frac{1}{2}\delta)$ which would mean that
  $\bigcup_i B(X_i,\frac{1}{2}\delta)$ is disconnected. Since we have $\log m$
  (Bernoulli) trials the probability for at least one to succeed is
  $1-\epsilon$. This proves the lemma and the theorem.
\end{proof}

\section{Future Work}
We believe it would be interesting to characterize the fractal sets
that naturally induce a semi uniform probability space. Also, in terms
of asymptotic analysis, it would be interesting to investigate other
geometric properties of these spaces, such as minimal matching, TSP
and power sum of MST edges.

\subsection*{Acknowledgements}
GK's work was supported by the National Science Foundation under
agreement no.~DMS-0111298. Any opinions, findings and conclusions or
recommendations expressed in this material are those of the authors
and do not necessarily reflect the views of the National Science Foundation.

\end{document}